\documentclass[11pt]{amsart}
\usepackage{amscd,amssymb,amsmath}
\newtheorem{theorem}{Theorem}[section]
\newtheorem{lemma}[theorem]{Lemma}

\theoremstyle{definition}

\theoremstyle{remark}
\newtheorem{remark}{Remark}[section]

\numberwithin{equation}{section}
\begin{document}
\title{Twisted conjugacy classes for polyfree groups}
\author{Alexander Fel'shtyn}
\address{Instytut Matematyki, Uniwersytet Szczecinski,
ul. Wielkopolska 15, 70-451 Szczecin, Poland}
\email {felshtyn@gmail.com}
\author{Daciberg Gon\c calves}
\address{Dept. de Matem\'atica - IME - USP, Caixa Postal 66.281 - CEP 05311-970,
S\~ao Paulo - SP, Brasil}
\email{dlgoncal@ime.usp.br}
\author{Peter Wong}
\address{Department of Mathematics, Bates College, Lewiston,
Maine 04240, U.S.A.}
\email{pwong@bates.edu}
\thanks{This work was initiated during the first and second authors' participation at the 5th International Siegen Topology Symposium, "Manifolds and their Mappings" Siegen, Germany, July 25 - 30, 2005}

\begin{abstract}
Let $G$ be a finitely generated polyfree  group. If $G$ has nonzero Euler characteristic then we show that $Aut(G)$ has a finite index subgroup in which every automorphism has infinite Reidemeister number. For certain $G$ of length 2, we show that the number of Reidemeister classes of every automorphism is infinite.
\end{abstract}
\date{\today}
\keywords{Reidemeister number, polyfree groups}
\subjclass[2000]{Primary: 20E45; Secondary: 55M20}
\maketitle

\newcommand{\af}{\alpha}
\newcommand{\et}{\eta}
\newcommand{\ga}{\gamma}
\newcommand{\ta}{\tau}
\newcommand{\ph}{\varphi}

\newcommand{\bt}{\beta}

\newcommand{\lb}{\lambda}

\newcommand{\wh}{\widehat}

\newcommand{\sg}{\sigma}

\newcommand{\om}{\omega}

\newcommand{\cH}{\mathcal H}

\newcommand{\cF}{\mathcal F}

\newcommand{\N}{\mathcal N}

\newcommand{\R}{\mathcal R}

\newcommand{\Ga}{\Gamma}

\newcommand{\cc}{\mathcal C}

\newcommand{\bea} {\begin{eqnarray*}}

\newcommand{\beq} {\begin{equation}}

\newcommand{\bey} {\begin{eqnarray}}

\newcommand{\eea} {\end{eqnarray*}}

\newcommand{\eeq} {\end{equation}}

\newcommand{\eey} {\end{eqnarray}}

\newcommand{\ovl}{\overline}

\newcommand{\vv}{\vspace{4mm}}

\newcommand{\lra}{\longrightarrow}

\bibliographystyle{plain}

\section{Introduction}

Let $\ph:G\to G$ be an endomorphism and let $R(\ph)$ denote the cardinality of the set of $\ph$-twisted conjugacy classes, i.e., the number of orbits of the action $\sigma \cdot \alpha \mapsto \sigma \alpha \ph(\sigma)^{-1}$ where $\sigma, \alpha \in G$, also known as the Reidemeister number of $\varphi$. In classical fixed point theory, $R(\varphi)$ plays an important role in estimating the Nielsen number of a selfmap $f$ on a compact manifold, where $\varphi$ is the homomorphism induced by $f$ on the fundamental group. In \cite{fel-hill}, A. Fel'shtyn and R. Hill conjectured that injective endomorphisms of a finitely generated group with exponential growth should have infinite Reidemeister number. For automorphisms, the conjecture holds for non-elementary Gromov hyperbolic groups (which include free groups of finite rank) as shown by G. Levitt and M. Lustig \cite{levitt-lustig} (see also \cite{fel:1}). While the conjecture does not hold in general (not even for automorphisms), it has generated recent interest in finding classes of groups for which the Reidemeister number is always infinite for automorphisms. We say that a group $G$ has the $R_{\infty}$ property if $R(\varphi)=\infty$ for every automorphism $\varphi \in Aut(G)$. It follows from \cite{levitt-lustig} that free groups of finite rank at least 2 have the $R_{\infty}$ property. Our primary interest of this paper is in the (in)finiteness of $R(\ph)$ for automorphisms $\ph$ of finitely generated polyfree groups.

\section{Polyfree groups}

Polyfree groups are natural extensions of free groups. In fact,
such groups arise naturally in topology as certain surface braid groups. More specifically the following groups are known to be polyfree:\\

a)   The pure braid group $P_n(S)$ for all $n>0$
 where $S$ is a compact surface with nonempty boundary. This follows from  E.
 Artin for $S=D$, the 2-disk, see \cite{ar}  or Proposition 10.11 of \cite{bu-zi},  and
 for the other surfaces see Theorem 1 of \cite{go-gu}.\\

b)  The full Artin braid group $B_n(D)$ for  $n=2,3$ and 4. For the case $n=2$, we have
$B_2(D)=\mathbb Z$. The cases $n=3$ and $4$ follow from \cite{go-li}. In fact, they show that $B_3(D)$ is isomorphic to $F_2(a,b)\rtimes \mathbb Z$ and $B_4(D)$ is isomorphic to
$((F_2(a,b)\rtimes F_2(u,v))\rtimes \mathbb Z)$ where $F_2(x,y)$ denotes the free group generated by $x$ and $y$.

\medskip

The class of virtually polyfree groups (meaning that groups which contains a subgroup of finite index which is polyfree) is interesting as it contains the braid groups of compact surfaces with nonempty boundary.

Recall that a group $G$ is (finitely generated) polyfree if it admits a subnormal series
$\frak R: 1= N_0 \triangleleft N_1 \triangleleft ....\triangleleft N_{k-1}\triangleleft N_k=G$ whose factors $F_i=N_i/N_{i-1}$ are free groups of finite rank $r_i$ and $\frak R$ is called a polyfree series of $G$. The length of the polyfree group $G$ is $k$, i.e., the number of nontrivial subgroups of a subnormal series and the Euler characteristic, denoted by $c$, is the product of the values $(r_i-1)$ for $i$ ranging from $1$ to $k$.

According to \cite{Me}, both $k$ and $c$ only depend on $G$ and are independent of the series $\frak R$.
The main result of \cite{Me} is:

\begin{theorem}$($\cite{Me}$)$ A polyfree group $G$ of length $k$ and Euler characteristic
$c\ne 0$ has only a finite number $N$ of distinct polyfree series and
$$ N\leq (c+1)^{(k-1)c+k^2-1}.$$
\end{theorem}

First let us recall the following well-known fact about the Reidemeister number of an endomorphism of a short exact sequence of groups.

\begin{lemma}\label{R-facts}
Consider the following commutative diagram
\begin{equation*}\label{general-Reid}
\begin{CD}
    1 @>>> A    @>>>  B @>>>    C @>>> 1 \\
    @.     @V{\ph'}VV  @V{\ph}VV   @V{\ovl \ph}VV @.\\
    1 @>>> A    @>>>  B @>>>    C @>>> 1
 \end{CD}
\end{equation*}
where the rows are short exact sequences of groups and the vertical arrows are group endomorphisms.
If $R(\ovl \ph)=\infty$ then $R(\ph)=\infty$. If $Fix \ovl{\ph}=1$ and $R(\ph')=\infty$ then $R(\ph)=\infty$.
\end{lemma}

Now we will derive the first consequence of the above result for our problem.

\begin{theorem}\label{general} Let the Euler characteristic of $G$ be different from zero. Then there is a subgroup of  finite index in $Aut(G)$ such that the Reidemeister number of every element of that subgroup is infinite.
\end{theorem}

\begin{proof} Let us consider the set of all subnormal series. The group $Aut(G)$ acts on this set. For any given $\varphi\in Aut(G)$ and a subnormal series $\frak R: 1= N_0 \triangleleft N_1 \triangleleft ....\triangleleft N_{k-1}\triangleleft N_k=G$, the image $\varphi(\frak R): 1= \varphi(N_0) \triangleleft \varphi(N_1) \triangleleft ....\triangleleft \varphi(N_{k-1})\triangleleft \varphi(N_k)=G$  is again a subnormal series. The isotropy of this action is a subgroup which is of finite index in $Aut(G)$. But every element of this subgroup provides an automorphism which preserves the subnormal series. In particular, every such element $\eta$ will give rise to an automorphism
of the short exact sequence $1\to N_{k-1} \to G\to G/N_{k-1} \to 1$. Since $G/N_{k-1}$ is a free group of finite rank greater than 1 because the Euler characteristic is non-zero, the Reidemeister number of the induced automorphism $\overline \eta$ on $G/N_{k-1}$ is infinite. It follows from Lemma \ref{R-facts} that $R(\eta)$ is infinite.
\end{proof}

In some sense Theorem \ref{general} is weak, since in fact we have not settled the $R_{\infty}$ question for the majority of the automorphisms. But perhaps more can be explored of the results of \cite{Me}. For the rest of the paper, we investigate the finiteness of the Reidemeister number of automorphisms of polyfree groups of length 2 without any assumption on the Euler characteristic.

\section{Polyfree groups of length 2}

The  simplest group in this family is  $\mathbb Z\times \mathbb Z$.
It is the only abelian group in this family  and  it is known  that there are  automorphisms such that the number of Reidemeister classes is finite.

Another group in this family is $\mathbb Z\rtimes \mathbb Z$ where the action is non-trivial.
This group is the fundamental group of the Klein bottle and it follows from \cite{daci-pet} or \cite{fel-daci} that this group has the $R_{\infty}$ property.

Given a finitely generated polyfree group $G$ of length 2, $G$ admits the following short exact sequence
$$
1\to F_r \to G\to F_s \to 1
$$
where $F_n$ denotes a free group of rank $n$. Since $F_s$ is free, this sequence splits so $G$ has the form $F_r \rtimes F_s$. To investigate the $R_{\infty}$ property for such groups, we proceed in several steps by considering the following cases:\\
Case 1. $F_r\times F_s$ where $\max\{r,s\}\geq 2$.\\
Case 2. $\mathbb Z\rtimes_{\theta} F_s$, $s\geq 2$.\\
Case 3. $F_2\rtimes_{\theta} \mathbb Z$.\\
Case 4. $F_r\rtimes_{\theta} \mathbb Z$, $r\geq 3$.\\
Case 5. $F_2\rtimes_{\theta} F_s$, $s\ge 2$.\\
Case 6. $F_r\rtimes_{\theta} F_s$ where  $\min\{r,s\}\geq 2$.\\

 Let $a\in F_n$ and denote by $Z_n(a)$ the centralizer of $a$ in the group $F_n$. It is well known that if $a$ is the trivial element then $Z_n(a)=F_n$; otherwise $Z_n(a)$ contains $a$ and is isomorphic to the infinite cyclic group. The group $Z_n(a)$ is equal to the  maximal  cyclic subgroup which  contains $a$.

Given an element of $F_r\times F_s$, it has the form $(a,b)$ where $a\in F_r$ and $b\in F_s$. We have:

\begin{lemma}\label{cent-1}
The centralizer of $(a,b)$ is:\\
{\rm I)} The group $F_r\times F_s$ if $a=b=1$.\\
{\rm II)} The subgroup $Z_r(a)\times F_s$ if $a\ne 1$ and $b=1$.\\
{\rm III)} The subgroup $F_r\times Z_s(b)$ if $a=1$ and $b\ne 1$.\\
{\rm IV)} The subgroup $Z_r(a)\times Z_s(b)$.
\end{lemma}

\begin{proof} The proof is straightforward and we leave it to the reader.
\end{proof}

\begin{theorem}\label{thm3.2} For any $r,s$ with $\max\{r,s\}\ge 2$, $F_r\times F_s$ has the $R_{\infty}$ property.
\end{theorem}
\begin{proof} Without loss of generality we may assume that $r\leq s$. First let $r < s$. Then for any automorphism $\varphi$, we have $\varphi(F_r)\subset F_r$ and $\varphi(F_s)\subset F_s$. To see that, we suppose that $\varphi(x)=(a,b)$ where $x\in F_r$ and $b\ne 1$. By Lemma \ref{cent-1}, the centralizer of $(a,b)$  is either $\mathbb Z\times \mathbb Z$ or $F_r\times \mathbb Z$. But the centralizer of $x$ is $\mathbb Z\times F_s$. Note that an automorphism maps the centralizer of an element isomorphically to the centralizer of the image of that element. So we get a contradiction
since neither $\mathbb Z\times \mathbb Z$ nor   $F_r\times \mathbb Z$ can be isomorphic to $\mathbb Z\times F_s$ $(r<s)$. We conclude that the automorphism $\varphi$ preserves the summands. Therefore the number of Reidemeister classes is infinite because it is infinite in at least one of the summands that is a free group of rank greater than 1.  Now suppose
$r=s$. Then either 1) $\varphi(F_r)\subset F_r$ and $\varphi(F_s)\subset F_s$, or 2)      $\varphi(F_r)\subset F_s$ and $\varphi(F_s)\subset F_r$.  As before, let $x\in F_r$ and suppose that $\varphi(x)=(a,b)$. By Lemma \ref{cent-1},  the centralizer of $(a,b)$  is either $\mathbb Z\times \mathbb Z$, $\mathbb Z\times F_s$, or   $F_r\times \mathbb Z$. But the centralizer of $x$ is $\mathbb Z\times F_s$. So it is not possible to have the centralizer of $(a,b)$  be equal to  $\mathbb Z\times \mathbb Z$ which means that $a$ and $b$ cannot be both non-trivial. Thus, either 1) or 2) holds.
In the first case, clearly the result follows as in the first part of the proof. In the second case, let us identify the Reidemeister classes of elements of the form $(1,b)$ as follows. Call $\varphi_1:F_r \to F_s$ and $\varphi_2:F_s \to F_r$. Let $\varphi_3=\varphi_1\circ\varphi_2$. The Reidemeister classes of $\varphi$ are in 1-1 correspondence with the Reidemeister  classes of the automorphism $\varphi_3: F_s \to F_s$; the number of which is infinite and the result follows.
\end{proof}

Now we consider groups of the form $\mathbb Z \rtimes_{\theta} F_s$, $s\ge 2$ as in Case 2. We will assume that the action $\theta:F_s \to Aut(\mathbb Z)=\mathbb Z_2$ is non-trivial otherwise $\mathbb Z \rtimes_{\theta} F_s =\mathbb Z \times F_s$ as in Case 1. Let $F_s'$ be the kernel which is a free group of rank $2s-1$.

\begin{lemma} \label{cent-2}
Let $(a,b)\in \mathbb Z \rtimes_{\theta} F_s$. The centralizer of $(a,b)$ is:\\
{\rm I)} The group $\mathbb Z\rtimes F_s$ if $a=b=1$.\\
{\rm II)} The subgroup $\mathbb Z\rtimes F_s'$ if $a\ne 1$ and $b=1$.\\
{\rm III)} The subgroup $\mathbb Z\times Z_s(b)$ if $a=1$ and $b \in F_s' $ and $ b \ne 1$.\\
{\rm IV)} The subgroup $Z_s(b)$ if $a=1$ and $b \in F_s-F_s'$.\\
{\rm V)} The subgroup $\mathbb Z\times Z_s(b)$  if
$a\ne 1 \ne b$ and   $Z_s(b)\subset F_s'$ or the subgroup $\mathbb Z\times (Z_s(b))^2$  if
$a\ne 1 \ne b  $ and   $Z_s(b)$ is not a subset of $ F_s'$.
\end{lemma}

\begin{theorem} For $s\ge 2$, the group $\mathbb Z\rtimes F_s$ has the $R_{\infty}$ property.
\end{theorem}
\begin{proof} It follows from Lemma \ref{cent-2} that the image of an element of
$\mathbb Z$ is sent by the automorphism to an element in $\mathbb Z$. So we obtain a homomorphism of short exact sequence
$$
0\to \mathbb Z \to \mathbb Z\rtimes F_s \to F_s \to 1.
$$
Since the quotient $F_s$ is free of rank greater than 1, the result follows from Lemma \ref{R-facts}.
\end{proof}

\section{Mapping tori of free groups of rank 2}

We now study the $R_{\infty}$ property for the family of groups of the form $F_r \rtimes \mathbb Z$ as in Cases 3 and 4. Let  $\phi:\mathbb Z \to Aut(F_r)$ be the action.

For an arbitrary $r>1$, we have:

\begin{lemma}\label{lem3} If $\phi(t)$, where $t$ is a generator of $\mathbb Z$,  is conjugation by an element $w\in F_r$, then $R(\varphi)$ is infinite for any automorphism $\varphi$.
\end{lemma}

\begin{proof} Consider the maximal cyclic subgroup of $F_r$ which contains $w$. Then it can be shown that $F_r\rtimes_{\phi} \mathbb Z$ has a center generated by $wt^{-1}$. Since the center is characteristic. We have an induced map in the  quotient of $F_r\rtimes_{\phi} \mathbb Z$ by its center. But this quotient is isomorphic to  a free group on $r$ generators. To see this, consider the composition $F_r \to F_r\rtimes \mathbb Z \to F_r\rtimes \mathbb Z/\langle wt^{-1}\rangle$. This map is clearly surjective. But it is also injective since $\langle wt^{-1}\rangle\cap F_r=\{1\}$.
\end{proof}

\begin{remark}\label{rmk1} We observe that the  condition above is equivalent to saying that the homomorphism $(\phi(t))^{ab}$ is the identity, if $r=2$. This follows from a result of W. Magnus, see \cite{kms}, section 3.5 Corollary N4,  p. 169. By \cite[Lemma 2.1]{bo-ma-ven}, $F_r \rtimes \mathbb Z \cong F_r \times \mathbb Z$ so Lemma \ref{lem3} follows from Theorem \ref{thm3.2}
\end{remark}

\begin{lemma}\label{lem4} If $R({\phi(t)}^{ab})<\infty$, then
$R(\varphi)=\infty$ for any automorphism $\varphi$.
\end{lemma}

\begin{proof} The commutator subgroup of $F_r\rtimes \mathbb Z$ is in fact a subgroup of $F_r$ which contains the commutator subgroup $[F_r,F_r]$. It is not difficult to see that the abelianization  is the sum of $\mathbb Z$ with the cokernel $Coker(1-{\phi(t)}^{ab})=\mathbb Z^r/(1-{\phi(t)}^{ab})(\mathbb Z^r)$. Since $R({\phi(t)}^{ab})<\infty$ and $R({\phi(t)}^{ab})=|Coker(1-{\phi(t)}^{ab})|$, it follows that the commutator subgroup of $F_r$ has finite quotient. This implies that it is a free group of finite rank greater than 1. Now using  the short exact sequence $$1\to [G,G] \to G \to G^{ab}\to 0$$
where $G=F_r\rtimes \mathbb Z$, the result follows from Lemma \ref{R-facts}.
\end{proof}

Now we focus on Case 3. We will make use of the following classification of the groups $F_2\rtimes_{\phi}\mathbb Z$ given by Theorem 1.1 of \cite{bo-ma-ven}

\begin{theorem}\label{bmv} $($\cite{bo-ma-ven}$)$ Let $F_2=\langle a,b\rangle$ be a free group of rank 2, let $\phi\in Aut(F_2)$, and consider the mapping torus $M_{\phi}=F_2\rtimes_{\phi}\mathbb Z.$ Let $\phi^{ab}\in GL_2(\mathbb Z)$ be the map induced by $\phi$ on $F_2^{ab}\cong \mathbb Z^2$ (written in row form with respect to $\{a,b\}$).\\

\item[i)]If $\phi^{ab}=I_2$, then $Out(M_{\phi})\cong (\mathbb Z^2\rtimes C_2)\rtimes GL_2(\mathbb Z)$, where $C_2$ is the cyclic group of order 2 acting on $\mathbb Z^2$ by sending $u \to -u$, $u\in \mathbb Z^2$, and where $GL_2(\mathbb Z)$ acts trivially on $C_2$ and naturally on $\mathbb Z^2$ (thinking vectors as columns there). \\

\item[ii)]  If $\phi^{ab}=-I_2$, then $Out(M_{\phi})\cong PGL_2(\mathbb Z)\times C_2$.\\

\item[iii)]  If $\phi^{ab}\ne -I_2$ and does not have 1 as an eigenvalue then  $Out(M_{\phi})$ is finite.\\

\item[iv)]  If $\phi^{ab}$ is conjugate to $\left(\begin{smallmatrix}
                                                      1 & k \\
                                                      0 & -1
                                                  \end{smallmatrix}\right)$ for some integer $k$, then $Out(M_{\phi})$ has an infinite cyclic subgroup of finite index.

\item[v)]  If $\phi^{ab}$ is conjugate to
$\left(\begin{smallmatrix}
        1 & k \\
        0 & 1
     \end{smallmatrix}\right)$ for some integer $k\ne 0$, then $Out(M_{\phi})$ has an infinite cyclic subgroup of finite index.\\

\noindent Furthermore, for every $\phi\in Aut(F_2)$, $\phi^{ab}$ fits into exactly one of the above cases.
\end{theorem}

We now determine the $R_{\infty}$ property for these mapping tori.

\begin{theorem}\label{polyfree2} The groups  $F_2\rtimes_{\phi} \mathbb Z$ have the $R_{\infty}$ property for any action $\phi$.
\end{theorem}
\begin{proof} According to the classification Theorem \ref{bmv}, we have 5 cases: i), ii), iii), iv) and v). \\
Case i) follows from Remark \ref{rmk1} that $\phi$ is conjugation together with Lemma \ref{lem3}. \\
Cases ii) and iii) follow from Lemma \ref{lem4}.\\
In Case iv), the homomorphism $\phi^{ab}$ has finite order. In fact
 $(\phi^{ab})^2=id$. This implies that $\phi^{\ell}$ induces the identity for some integer $\ell$. So by \cite{kms} section 3.5 Corollary N4 (p. 169),
  $\phi^{\ell}$
is conjugation by some word $w$. This means that the automorphism $w^{-1}\phi w$ has finite order in $Out(F_2)$. By a G(eneralized) B(aumslag)-S(olitar) group, we mean the fundamental group of a graph of $\mathbb Z$'s. Let $\Gamma$ be the mapping torus obtained using the automorphism $w^{-1}\phi w$. Since $w^{-1}\phi w$ has finite order in $Out(F_2)$, it follows that $F_2 \times \mathbb Z$ is a finite index subgroup in $\Gamma$ so that $\Gamma$ is quasi-isometric to $F_2 \times \mathbb Z$, a non-elementary GBS group, and hence, by \cite{TW}, has property $R_{\infty}$. Alternatively, using the classification of graphs of $\mathbb Z$'s \cite{whyte}, $\Gamma$ is in fact a non-elementary GBS group.
It follows from \cite{lev} that $\Gamma$ has the $R_{\infty}$ property. But the two mapping tori $\Gamma$ and $F_2\rtimes_{\phi} \mathbb Z$ are isomorphic. So the result follows.\\

\medskip

In case v), the matrix of  $\phi^{ab}$ is of the form $\left(\begin{smallmatrix}
                                                           1 & k \\
                                                           0 & 1
                                                       \end{smallmatrix}\right)$
 for some nonzero integer $k$.
 By composing with an appropriate inner automorphism in $F_2$, we may assume that
the automorphism $\phi$ is given by $\phi(x)=xy^k$ and $\phi(y)=y$. Now, we have
$$M_{\phi}=\langle x,y,t|t^{-1}xt=xy^k, t^{-1}yt=y\rangle.
$$
By the proof of  Proposition 4.4 in \cite{bo-ma-ven},
every automorphism $\theta\in Aut(M_{\phi})$ is of one of the
following
 forms: $a')$  $\Psi^{-m}\gamma_g^{-1}\Xi^i$,
$b')$ $\Psi^{-m}\gamma_g^{-1}\Omega\Xi^i$, $c')$ $\Psi^{-m}\gamma_g^{-1}\Delta\Xi^i$,
$d')$  $\Psi^{-m}\gamma_g^{-1}\Omega\Delta\Xi^i$ for $m$ an integer and $0\leq
i\leq |k|-1$. Here, $\Psi, \Omega, \Delta$, and $\Xi$ are given by
\begin{equation*}
\begin{aligned}
  &\Psi        &\quad        &\Omega         &\quad        &\Delta           &\quad       &\Xi \\
x &\mapsto tx  &\quad      x &\mapsto x      &\quad      x &\mapsto x^{-1}   &\quad     x &\mapsto xy  \\
y &\mapsto y   &\quad      y &\mapsto y^{-1} &\quad      y &\mapsto y^{-1}   &\quad     y &\mapsto y   \\
t &\mapsto t   &\quad      t &\mapsto t^{-1} &\quad      t &\mapsto ty^{-k}  &\quad     t &\mapsto t
\end{aligned}
\end{equation*}
and $\gamma_g$ is conjugation by $g$. Note that the notation used in \cite{bo-ma-ven} is composition on the right.

Since the induced automorphism in an abelian quotient remains the same if we
change
an automorphism by conjugation, we will reduce our problem to considering the
automorphisms
$a)$  $\Psi^{-m}\Xi^i$,   $b)$ $\Psi^{-m}\Omega\Xi^i$, $c)$
$\Psi^{-m}\Delta\Xi^i$,
$d)$  $\Psi^{-m}\Omega\Delta\Xi^i$ for $m$ an integer and $0\leq i\leq |k|-1$.

Equivalently, the automorphisms in question are given as follows:
$a)$ $x \mapsto t^{-m}xy^i, y \mapsto y,  t \mapsto t$;
$b)$ $x \mapsto t^{-m}xy^{-i}, y \mapsto y^{-1}, t \mapsto t^{-1}$;
$c)$ $x \mapsto x^{-1}t^my^{-i}, y \mapsto y^{-1}, t \mapsto ty^{-k}$; and
$d)$  $x \mapsto x^{-1}t^my^i, y \mapsto y, t \mapsto t^{-1}y^k$.

By sending $x \mapsto \ovl{x}, y \mapsto \ovl 1, t \mapsto \ovl{t}$, we have a short exact sequence
$$1\to H \to G \to \mathbb Z\times\mathbb Z \to 0$$
where $G=F_2\rtimes_{\phi} \mathbb Z$ and $H$ is the normal closure of $y$ in $G$. Equivalently, this projection is given by $w \mapsto (w_x,w_t)$ where $w\in G$ is a word and $w_x$ and $w_t$ denote the $x$-exponent sum and the $t$-exponent sum of $w$ respectively.
It is not difficult to see that the automorphisms in the four cases induce automorphisms on the quotient and the induced maps $\ovl{\theta}$ are
given by
the following matrices (using the convention of \cite{bo-ma-ven}):
\begin{displaymath}
a) \left(\begin{smallmatrix}
             1 & -m \\
             0 & 1
             \end{smallmatrix}\right) \quad b)\left(\begin{smallmatrix}
                                               1 & -m \\
                                               0 & -1
                                               \end{smallmatrix}\right) \quad c)\left(\begin{smallmatrix}
                                                                                   -1 & m \\
                                                                                   0 & 1
                                                                                \end{smallmatrix}\right) \quad d)\left(\begin{smallmatrix}
                                                                                           -1 & m \\
                                                                                           0 & -1
                                                                                          \end{smallmatrix}\right).
\end{displaymath}
In cases $a), b), c)$, we see that the Reidemeister number $R(\ovl{\theta})=\infty$ since $\det (I-\ovl{\theta})=0$.

For case $d)$ however, $R(\ovl{\theta})<\infty$ so we need to use a different argument as follows. Since $G=F_2(x,y) \rtimes \langle t\rangle$, the subgroup $H$ is a subgroup of $F_2(x,y)$, so it is a free group. By the Reidemeister-Schreier method, $H$ is generated by the set $\{x^iyx^{-i} \ | i\in \mathbb Z \}$. The projection of this set
onto the abelianization of $H$ yields a generating set for the free abelian group $H^{ab}$. Consider the  induced automorphism $\eta=\theta|_H$ on $H$. Then the induced homomorphism $\eta^{ab}$ on $H^{ab}$ has infinite Reidemeister number.
To see that, we recall that the presentation of $M_{\phi}$ is $\langle x,y,t |  t^{-1}xt=xy^k, t^{-1}yt=y \rangle$,
and the automorphism $\theta$ is given by $x \mapsto x^{-1}t^my^i, y \mapsto y, t \mapsto t^{-1}y^k$. Let $H_0$ be the subgroup of $H$ generated by the element $y$ and $\ovl{H_0}$ be the image of $H_0$ in $H^{ab}$.
Now $\eta$ induces the identity on $\ovl{H_0}$ and thus the number of Reidemeister classes of
the automorphism $\eta$ restricted to $H_0$ is infinite. It follows that $R(\eta^{ab})=\infty$ and so the Reidemeister number $R(\eta)$ is infinite. Note that $R(\ovl{\theta})<\infty$ since $\det (I-\ovl{\theta})\ne 0$. It follows that $Fix \ovl{\theta}$ is trivial and hence $R(\theta)=\infty$ by Lemma \ref{R-facts}. To complete the proof, we need to show that the automorphism $\Theta=\Psi^{-m}\gamma_g^{-1}\Omega\Delta\Xi^i$ as in case $d')$ also has infinite Reidemeister number. Since conjugation $\gamma_g$ does not change the exponent sums of $x$ and $t$, the subgroup $H$ is invariant under $\Theta$. Moreover, $\Theta$ and $\theta$ induce the same automorphism $\eta^{ab}$ on the abelianization $H^{ab}$. It follows that $R(\Theta)$ is also infinite.
\end{proof}

\section{Conclusion}

For the remaining cases, namely Cases 4: $F_r\rtimes \mathbb Z$ for $r> 2$,   5: $F_2\rtimes F_s$ for $s\geq 2$ and      6: $F_r\rtimes F_s$ for $\min\{r,s\}\geq 2$ we are able to decide only the subcases of case 4 for which the hypotheses of Lemma \ref{lem3} or Lemma \ref{lem4} are satisfied. So we will keep these Cases 4,5 and 6 as open questions.

\noindent
{\bf Acknowledgments.} The first author is grateful to Zoran \v Sunic for very helpful discussion and for bringing \cite{Me} to his attention. The second author would like to thank A. Martino for many helpful discussion about automorphisms of free-by-cyclic groups during a workshop in Lewiston, Maine, April 2006 and the AMS-PTM joint meeting in Warsaw, August 2007. The third author would like to thank E. Ventura for pointing out an error in an earlier proof of Case iv) of Theorem \ref{polyfree2}.

\end{document}